\documentclass[12pt,letterpaper]{amsart}
\usepackage{amsmath, amsthm, amssymb}
\usepackage[tmargin=1.4in,bmargin=1.2in,rmargin=1.4in,lmargin=1.4in]{geometry}
\usepackage[breaklinks=true]{hyperref}
\usepackage{soul}
\usepackage{changes}
\usepackage{enumitem}

\theoremstyle{plain}
\newtheorem{theorem}{Theorem}[section]
\newtheorem{lemma}{Lemma}[section]
\newtheorem{prop}{Proposition}[section]
\newtheorem{cor}{Corollary}[section]

\theoremstyle{definition}
\newtheorem{defin}{Definition}[section]

\DeclareMathOperator{\tr}{tr}
\DeclareMathOperator{\End}{End}
\DeclareMathOperator{\Ima}{Im}
\DeclareMathOperator{\Id}{Id}

\DeclareMathOperator{\rank}{rank}
\begin{document}

\title{Representations of the Weyl algebra modulo $p^2$}
\author{Ruben Mamani-Velasco, Akaki Tikaradze}
\email{RubenLimbert.MamaniVelasco@rockets.utoledo.edu, Akaki.Tikaradze@utoledo.edu}
\address{University of Toledo, Department of Mathematics \& Statistics, 
Toledo, OH 43606, USA}

\begin{abstract}
Let $\mathbf{k}$ be an algebraically closed field of characteristic $p>0$. Let $R=W_2(\mathbf{k})$ be the ring of Witt vectors of length $2$.
In this note, we are concerned with the study of representations of the first Weyl algebra that are finite-rank free
$R$-modules ($R$-lattices). 
Our main result shows that given an $R[x]$-lattice $V$, it admits a left module structure over
the Weyl algebra (equivalently, an $R$-linear connection) if and only if all of its indecomposable summands have 
the form $R[x]\otimes_{R[x^p]}N$, where $N$ is an indecomposable $R[x^p]$-lattice of rank divisible by $p$.  
 We also classify, up to isomorphism, all lattices over
the $p$-Weyl algebra (defined by $[y,x]=p$) for which the action of $x$ modulo $p$ is regular nilpotent.
\end{abstract}

\maketitle

\section{Introduction and main results}
Let $R$ be a unital commutative ring. Throughout, given an $R$-algebra $A$, by an $A$-lattice we mean a left $A$-module which is a finite-rank free $R$-module.

Let $A_1(R)=R\langle x,y\rangle/(yx-xy-1)$ denote the first Weyl algebra over $R$. The study of modules over the Weyl algebra has been a topic of considerable interest in representation theory and noncommutative algebra for some time now. While the cases in which the base ring $R$ is an algebraically closed field (both in characteristic 0 and positive characteristic) are essentially
well-understood, much less is known when $\mathbb{Z}/p^n\mathbb{Z}\subset R$, where $p$ is a prime and $n>1$. 

From now on,  we fix an algebraically closed field $\mathbf{k}$ of characteristic $p>0$ and set $R=W_n(\mathbf{k})$  to be the ring of Witt vectors of length $n\geq 1$.  In fact, most of our results are in the case $n=2$. 

We are primarily interested in studying $A_1(R)$-modules which are finite free $R$-modules; we call such modules
$A_1(R)$-lattices.  Specifically, we are interested in determining if
a given $R[x]$-lattice admits an $A_1(R)$-module structure (equivalently an $R$-linear connection).
 Since the reduction modulo $p$ of such a module is an $A_1(\mathbf{k})$-module finite-dimensional over $\mathbf{k}$, 
it is important to understand which $A_1(\mathbf{k})$-modules
can be lifted to an $A_1(R)$-lattice. The answer is provided by the following theorem, which is the main result of this paper.

We first introduce a notion of a quasi-Kn\"{o}rr lattice, directly motivated by
the notion of a virtually irreducible lattice introduced by R. Kn\"{o}rr  in his work on Brauer's height-zero conjecture \cite{K} (see also \cite{G}).

\begin{defin}
Let $R$ be a local Artinian ring. Let $A$ be an $R$-algebra. Let $M$ be an $A$-lattice. 
Then we call $M$ a quasi-Kn\"{o}rr $A$-lattice if it is indecomposable,
and for any $z\in \End_A(M)$  we have   $\tr(z)\in (\rank_{R}(M))$. 
    
\end{defin}

\begin{theorem}\label{Main-Weyl}
 Let $M$ be an $R[x]$-lattice. Then it admits a left $A_1(R)$-module structure if and only if its indecomposable summands are quasi-Kn\"{o}rr
 lattices of rank divisible by $p^n$. If $n=2$,  then $M$ admits a left $A_{1}(R)$-module structure if and only if all its indecomposable summands
 are of the form $R[x]\otimes_{R[x^p]}N$,  where $N$ is an indecomposable $R[x^p]$-lattice of rank divisible by $p$. 
 Let $V$ be a finite-dimensional $A_1(\mathbf{k})$-module. Then there exists an $A_{1}(R)$-lattice $M$ (still $R=W_2(\mathbf{k})$) such that $M/pM\cong V$ as
 $A_1(\mathbf{k})$-modules if and only if all indecomposable summands of $V$ have $\mathbf{k}$-dimension divisible by $p^2$.

\end{theorem}

An important ingredient for the proof of the above result is the equivalence between the category of $A_1(R)$-lattices and the category
of lattices over the $p$-Weyl algebra, defined as the algebra over $R$ generated by $y,x$ subject to the relation $[y,x]=p$.
We denote the $p$-Weyl algebra over $R$ by $A_{1,p}(R)$.

 Finally, 
we give the full classification of all non-isomorphic $p$-Weyl module structures on an $R[x]$-lattice $M$
under the assumption that the action of $x$ on $M/pM$ is regular nilpotent, Theorem \ref{thm:Regular}.

We also note that a closely related question about describing two-sided ideals of $A_1(R)$  which form an $R$-lattice, has been considered in
\cite{T}. In particular, we showed that all such ideals are generated by central elements.

\section{Proof of main theorem}

We first recall some well-known facts about finite-dimensional representations of $A_1(\mathbf{k})$. The Weyl algebra
is a rank $p^2$ Azumaya algebra over its center $\mathbf{k}[x^p, y^p]$.  Thus all its irreducible representations are $p$-dimensional
and uniquely determined by their central characters. The category of finite dimensional $A_1(\mathbf{k})$-modules breaks up as the direct sum of subcategories of modules supported on the closed points of the spectrum of the center. Thus, it is enough to describe
the category of finitely generated $A_1(\mathbf{k})/\big((x^{p}-\lambda)^m,(y^{p}-\mu)^m\big)$-modules for $m\geq 1$ and $\lambda, \mu\in \mathbf{k}$.

For this, one needs to describe a splitting module $L$ of $A_1(\mathbf{k})/\big((x^{p}-\lambda)^m, (y^{p}-\mu)^m\big)$ (unique up to an isomorphism). Then
this category is equivalent to the category of finitely generated $\mathbf{k}[x^p, y^p]/\big((x^p-\lambda)^m,(y^p-\mu)^m\big)=Z$-modules,

and the equivalence is given by
$$M\to L\otimes_Z M,\quad M \in Z\text{-mod}.$$ 

We also note that the descriptions of the splitting modules mentioned above are closely related to the splitting module
of $A_1(\mathbf{k})/(y^p)$  (which is just the standard $A_1(\mathbf{k})$-module $\mathbf{k}[x])$.  This fact is a particular case

of the Cartier descent which we briefly recall next. Let $X$ be a smooth variety over $\mathbf{k}$. The ring $D(X)$ of crystalline differential operators on $X$
is an Azumaya algebra over the Frobenius twist of the cotangent bundle of $X$ \cite{BMR}.
Then $D(X)$ splits as an Azumaya algebra over the zero section of the cotangent bundle, which yields
 an equivalence of categories of $\mathcal{O}_X$-modules with a flat connection with vanishing
$p$-curvature and the category of modules over the Frobenius twist of $X$,
this is the general Cartier descent (for the Weyl algebra case, see, for example, \cite[Lemma~2.3]{Ts}). This result was significantly generalized further by

Ogus-Vologodsky \cite{OV} for a nilpotent $p$-curvature under some assumptions on $X$
(which hold for the local case we are interested in).

The situation is considerably more complicated when it comes to modules over the Weyl algebra over 
$\mathbb{Z}/p^n\mathbb{Z}$ for $n\geq 2$.  Indeed, the Weyl algebra $A_1(R)$ is no longer Azumaya over its center (which is still large).
A partial analog of the above results holds. Shiho \cite{S} showed that if $X$ is a smooth variety over $R$ admitting a lift of the Frobenius (which is the case for an affine $X$), then the category of coherent sheaves on $X$ equipped with a nilpotent flat connection is equivalent

to the category of coherent sheaves on the Frobenius twist of $X$ equipped with a nilpotent flat $p$-connection.

Our starting point is the following relation between modules over the Weyl and $p$-Weyl algebras, which in particular implies that the classification problems of $R$-lattices over these algebras are equivalent. In what follows $A_1(R)^f$-mod (respectively $A_{1,p}(R)^f$-mod) denotes the full subcategory
of $A_1(R)$-modules (respectively, $A_{1,p}(R)$-modules) that are finite as $R$-modules.
\begin{prop}\label{Equivalence}
 There exists an $R$-linear equivalence of categories 
 $$F:A_{1,p}(R)^f\text{-mod}\to A_{1}(R)^f\text{-mod},$$ such that $F(N)\cong R[x]_{\mathrm{Fr}}\otimes_{R[x]}N$ as $R[x]$-modules, where $\mathrm{Fr}:R[x]\to R[x]$ 

 is the $R$-algebra endomorphism determined by $\mathrm{Fr}(x)=x^p$.
\end{prop}
\begin{proof}
Recall that $x^{p^n}, y^{p^n}$ belong to the center of $A_1(R)$.  Similarly, $x^{p^{n-1}}, y^{p^{n-1}}$  are central in $A_{1, p}(R)$. 

Given $\alpha\in \mathbf{k}$, we denote by $A_1(R)^{f,\alpha}$-mod the full subcategory of $A_1(R)^f$-modules
on which $y^{p^n}-\tilde{\alpha}^{p^n}$ acts nilpotently, where $\tilde{\alpha}\in R$ is any lift of $\alpha$. 
It is clear that this definition is independent of a choice of a lift $\tilde{\alpha}$. 
We define the category $A_{1,p}(R)^{f,\alpha}$-mod as the full subcategory
of finite $A_{1,p}(R)$-modules on which $y-\tilde{\alpha}$ acts nilpotently. This definition is also independent of the choice of the lift $\tilde{\alpha}$. It is immediate that

the category $A_1(R)^f\text{-mod}$ breaks up as the  direct sum of $A_1(R)^{f,\alpha}\text{-mod}$ as $\alpha\in\mathbf{k}$ varies. A similar decomposition holds for
$A_{1,p}(R)^f\text{-mod}$. Now, the existence of the desired equivalence between $A_1(R)^{f,0}\text{-mod}$ and $A_{1,p}(R)^{f,0}\text{-mod}$
follows from Shiho's equivalence discussed above. Finally, by translating $y$ by a lift $\tilde{\alpha}$ on both sides, we obtain an equivalence
$A_1(R)^{f,\alpha}\text{-mod}\to A_{1,p}(R)^{f,\alpha}\text{-mod}$. Taking the direct sum of these equivalences over $\alpha\in\mathbf{k}$, 

we are done.
\end{proof}

Let us start with a few preliminary lemmas which should be fairly standard. We include the proofs for the reader's convenience.

\begin{lemma}\label{Frobenius-indecomposable}
Let $N$ be an $R[x^p]$-lattice. Then $N$ is indecomposable if and only if $R[x]\otimes_{R[x^p]}N$ is an indecomposable $R[x]$-lattice.
    
\end{lemma}
\begin{proof}
 Let $N$ be an indecomposable $R[x^p]$-lattice. Put $A=\End_{R[x^p]}(N)$. Hence $A$ has no nontrivial idempotents. 
So $A/J(A)\cong\mathbf{k}$, where $J(A)$ denotes the Jacobson radical. Let $\alpha\in \mathbf{k}$  such that $x^p\mod J(A)=\alpha$. 
 Put $B=A\otimes_{R[x^p]}R[x]$,  then $B=\End_{R[x]}(R[x]\otimes_{R[x^p]}N)$.  Thus, we need to show that $B$ has no nontrivial idempotents, equivalently $B/J(B)=\mathbf{k}$. Choose $\beta\in R$ such that $\bar\beta^{\,p}=\alpha$, where $\bar\beta$ denotes the reduction of $\beta$ modulo $p$. Then $x-\beta\in J(B)$,
 so $B=R+J(B)$. So $B/J(B)=\mathbf{k}$  as desired.

\end{proof}

\begin{lemma}\label{lemma:trace}
Let $M$ be a finite free $R$-module. Let $c\in R$ and $x_i \in \End_R(M)$, $1\leq i\leq m$.
There exist $y_i \in \End_R(M)$ such that
\[
\sum_{i=1}^m[y_i, x_i] = c \cdot \Id_M,
\]
if and only if for any $z \in \End_{R}(M)$ commuting with all $x_i$, we have $\tr(cz)=0$. 

\end{lemma}

\begin{proof}
At first, recall that the trace pairing  
$$\langle -, -\rangle: \End_R(M)\times \End_R(M)\to R: \langle A, B\rangle=\tr(AB)$$ 
induces an $R$-module isomorphism 
$\End_R(M)\cong \End_R(M)^*$, here $N^*=\operatorname{Hom}_R(N, R)$  for an $R$-module $N$. 

For an $R$-submodule $N\subset \End_R(M)$ , denote by $N^{\perp}$ its orthogonal complement under the trace pairing.
The key part of the proof is showing the double orthogonal complement property: for any $R$-submodule $N$ of $\End_R(M)$,  we have
$(N^{\perp})^\perp=N$. Indeed, we first remark that for any finitely generated $R$-module $N$,  we have
$$N\cong N^*,\qquad \operatorname{Ext}^1_R(N,R)=0.$$

Indeed, every finite $R$-module is a direct sum of cyclic modules $R/p^iR$, for which these assertions follow directly.
Consequently, the contravariant functor $N\mapsto N^*$ is exact on the category of finite $R$-modules.
Thus, by dualizing the exact sequence 
$$0\to N \to \End_R(M)\to Q\to 0,$$  and using the trace pairing, we
may identify $N^{\perp}$ with $Q^*$. Dualizing one more time, we get that
 $(N^{\perp})^{\perp}\cong N^{**}$. Thus  $(N^{\perp})^{\perp}$ and $N$ have equal lengths as $R$-modules, and since
 $N\subset (N^{\perp})^{\perp}$, we get the desired equality $N=(N^{\perp})^{\perp}$.

Let $z\in \End_R(M)$ commute with each $x_i$, $1\leq i\leq m$. Then we have
\[
\sum_{i=1}^m[zy_i,x_i]=cz.
\]
Taking the trace of both sides, we obtain the necessity of the condition.
Next, we show that the condition is sufficient.
Define the $R$-linear map
\[
\varphi : \prod_{i=1}^m\End_R(M) \longrightarrow \End_R(M), \qquad \varphi(y_1,\dots,y_m) = \sum_{i=1}^m[y_i, x_i].
\]
Thus, we need to show that $c \Id_M \in \Ima(\varphi)$.

A direct computation shows that the orthogonal complement of $\Ima(\varphi)$ is
\[
\Ima(\varphi)^\perp=\End_{R\langle x_1,\ldots,x_m\rangle}(M).
\]

By the double orthogonal complement result, $c \Id_M \in \Ima(\varphi)$ if and only if $c \Id_M$ is orthogonal to $\End_{R\langle x_1,...,x_m\rangle}(M)$. That is, for every $z\in \End_{R\langle x_1,...,x_m\rangle}(M)$,
\[
c \cdot \tr(z)=\tr(c \Id_M \cdot z) = 0.
\]
Therefore, there exist $y_1,...,y_m$ satisfying the desired equality if and only if $c \cdot \tr(z) = 0$ for all $z \in \End_{R\langle x_1,...,x_m\rangle}(M)$.
\end{proof}
\begin{lemma}\label{lem:commuting-lifts}
Let $M$ be an $R$-lattice. Put $V=M/pM$. 
Let $\bar x,\bar y\in\operatorname{End}_{\mathbf{k}}(V)$ be commuting endomorphisms. 
Then there exist commuting endomorphisms $x_0,y_0\in\End_R(M)$ whose reductions modulo $p$ are $\bar x$ and $\bar y$, respectively.
\end{lemma}

\begin{proof}
First, let us assume $\bar x$ is nilpotent. Let $d=\dim_{\mathbf{k}}V$ and choose a basis in which $\bar x$ has Jordan normal form. Thus
\[
\bar x = \sum_{i=1}^{d-1} c_i e_{i,i+1},\qquad c_i\in\{0,1\},
\]
where $e_{ij}\in \End_{\mathbf{k}}(V)$  are the standard elementary matrices $1\leq i, j\leq d$. 
Lift $\bar x$ to an element $x_0\in M_d(R)$ by taking the same matrix with entries $c_i$ (viewed in $R$). Then $x_0=\bar x\mod p$.

Write $\bar y=\sum_{i,j} a_{ij}e_{ij}$ with $a_{ij}\in\mathbf{k}$. The commuting condition $[\bar x,\bar y]=0$ is equivalent to the linear system
\[
c_r\, a_{r+1,s} = c_{s-1}\, a_{r,s-1}\qquad (1\le r,s\le d),
\]
with the convention $c_0=c_d=0$ and $a_{i,j}=0$ if $i$ or $j$ is out of range.

This system has a very special form: each equation says that each entry $a_{ij}$ is either $0$, or forced to be equal to another entry. Indeed,

\begin{itemize}
\item if $c_r=0$, the equation becomes $0=c_{s-1} a_{r,s-1}$, so either $c_{s-1}=0$ (no condition) or $a_{r,s-1}=0$;
\item if $c_r=1$, then $a_{r+1,s}=c_{s-1} a_{r,s-1}$, so either $a_{r+1,s}=0$ (if $c_{s-1}=0$) or $a_{r+1,s}=a_{r,s-1}$.
\end{itemize}
Thus, we choose lifts of $a_{ij}$  to $R$ satisfying the same system of equations. This defines a matrix $y_0\in M_d(R)$ with $y_0=\bar y\mod p$ and  $[x_0,y_0]=0$.

Now, we consider the general case in which $\bar x$ might not be nilpotent. As in the proof of Proposition \ref{Equivalence}, $M$ can be written as a direct sum of $R$-modules such that the reduction modulo $p$ of each summand is invariant under $\bar x$ and $\bar y$, and there exists a $\lambda\in\mathbf{k}$ such that $\bar x-\lambda$ acts nilpotently on it. Thus, we can find commuting lifts for each summand and we are done.

\end{proof}

\begin{lemma}\label{indecomposable-trace}
Let $A$ be an $R$-algebra, let $M$ be an indecomposable $A$-lattice. Suppose that $\rank_R(M)$  is divisible by $p$. 
Then for any $z\in \End_A(M)$,  we have $\tr(z)\in (p)$.

\end{lemma}
\begin{proof}
Suppose that there exists $z\in \End_A(M)$,  so that $\tr(z)\neq 0 \mod p$. Let $\bar{z}_{s}$  denote the semisimple part
of $\bar{z}=z\mod p\in \End_{\mathbf{k}}(\overline{M})$. Here $\overline{M}=M/pM$. Then $\tr(\bar{z})\neq 0$. 

It follows that $\bar{z}_{s}$  has an eigenspace with a nonzero eigenvalue of dimension not divisible by $p$. Let $e$ be the projection on that eigenspace of
$\overline{M}$.  Then $e\in \mathbf{k}[\bar{z}]$  and $e^2=e$.  Let $\hat{e}\in R[z]$  be a lift of $e$ which is an idempotent. Then $\tr(\hat{e})\neq 0 \mod p$
and $\hat{e}\in \End_A(M)$,  a contradiction.

\end{proof}

\begin{proof}[Proof of Theorem \ref{Main-Weyl}]
Let $M$ be an $R[x]$-lattice. Then by Lemma \ref{lemma:trace} $M$ admits a left $A_1(R)$-module structure if and only
if for all $z\in \End_{R[x]}(M)$  we have $\tr(z)=0$. It is immediate that this condition holds if and only if the similar trace condition holds for all direct summands of $M$.  Thus, without loss of generality, we may assume that $M$ is an indecomposable $R[x]$-lattice.
So, as $\tr(\Id_M)=\rank_R(M)$,  we conclude that $p^n$ divides $\rank_R(M)$  and $\tr(z)\in (\rank_R(M))$.  Thus, $M$ admits a structure of $A_1(R)$-module if and only if all indecomposable summands of $M$ are quasi-Kn\"{o}rr of rank divisible by $p^n$.

We now assume $n=2$. We claim that an indecomposable $R[x]$-lattice $M$ of $R$-rank divisible by $p^2$ admits an $A_1(R)$-module structure
if and only if $M\cong R[x]\otimes_{R[x^p]}N$, for some indecomposable $R[x^p]$-lattice $N$.   Using Proposition \ref{Equivalence}, this is equivalent to proving that an indecomposable 
$R[x]$-lattice  admits a structure of $p$-Weyl module if and only if its rank is a multiple of $p$. Using Lemma \ref{lemma:trace}, this is equivalent to showing that $\tr(z)\in (p)$ for all $z\in \End_{R[x]}(N)$. This follows immediately from Lemma \ref{indecomposable-trace}.

Finally, to show the lifting criterion for a finite-dimensional $A_1(\mathbf{k})$-module, it is sufficient to prove the analogous statement for
existence of liftings for $A_{1, p}(R)$-modules (by Proposition \ref{Equivalence}).
In what follows, $\mathbf{k}[\bar{x},\bar{y}]$ is a polynomial ring identified with $A_{1,p}(R)/pA_{1,p}(R)$, where $\bar{x}=x\bmod p$ and $\bar{y}=y\bmod p$. 
We need to show that given a finite-dimensional $\mathbf{k}[\bar{x}, \bar{y}]$-module $V$, 
there exists an $A_{1,p}(R)$-lattice $M$ and an isomorphism of $\mathbf{k}[\bar{x}, \bar{y}]$-modules
$M/pM\cong V$ if and only if all indecomposable summands of $V$ have dimensions divisible by $p$. 
 Choose a finite free $R$-module $M$ together with an identification $M/pM\cong V$ as $\mathbf{k}$-vector spaces.
 By Lemma \ref{lem:commuting-lifts}, there exist commuting $x_0,y_0\in\End_R(M)$ lifting $\bar{x},\bar{y}$, respectively. 
 Let $x, y\in \End_R(M)$  be any lifts of $\bar{x}, \bar{y}$  respectively.
 Therefore, there exist $\widetilde A,\widetilde B\in\End_R(M)$ such that $x=x_0+p\widetilde A$ and $y=y_0+p\widetilde B$. Let $A,B\in\End_{\mathbf{k}}(V)$ denote the reductions of $\widetilde A,\widetilde B$ modulo $p$. Since $p^2=0$ and $[x_0,y_0]=0$, the relation $[y,x]=p\Id_M$ holds if and only if
$$[\bar{x}, -B]+[\bar{y}, A]=\Id_V$$

By Lemma \ref{lemma:trace}, there exists $A, B$ satisfying the above equation if and only if for every $z\in\End_{\mathbf{k}[\bar{x},\bar{y}]}(V)$, 
we have $\tr(z)=0$. Now, using Lemma \ref{indecomposable-trace}, we are done.

\end{proof}
As a corollary of the proof above (and Lemma \ref{Frobenius-indecomposable}), we have obtained the following classification of quasi-Kn\"{o}rr lattices of rank divisible by $p^2$ (for $R=W_2(\mathbf{k})$),
which might be of independent interest.

\begin{cor}
Let $R=W_2(\mathbf{k})$. Let $V$ be an indecomposable $R[x]$-lattice of rank divisible by $p^2$.  Then $V$ is a quasi-Kn\"{o}rr $R[x]$-lattice
if and only if there exists an indecomposable $R[x^p]$-lattice $N$ of rank divisible by $p$ such that $V\cong R[x]\otimes_{R[x^p]}N$.

\end{cor}

\section{Regular nilpotent case}
In this section, we assume $R=W_2(\mathbf{k})$. 
The goal of this section is to classify all $A_{1,p}(R)$-lattices under the additional assumption that the action of $x$ modulo
$p$ is regular nilpotent. We first introduce a family of $A_{1,p}(R)$-lattices.

\begin{defin}
Let $m\geq 1$ and $f(x)\in \mathbf{k}[x], g(x)\in R[x]$ with $\deg f, \deg g\leq mp-1$. Then an $A_{1,p}(R)$-lattice $M_{f,g}$ of rank $mp$  is defined as follows. As an $R[x]$-module we put $M_{f, g}=R[x]/(x^{mp}-pf(x))$. The action of $y$ is given by $y=p\partial+g(x)$, where $\partial=\frac{d}{dx}$. Here $pf(x)$ means $p\widetilde f(x)$ for any coefficientwise lift $\widetilde f(x)\in R[x]$ of $f(x)$; this is independent of the choice since $p^2=0$.

\end{defin}

We also introduce the following subset of $\mathbf{k}[x]/(x^{mp})$, representing each element by its unique polynomial of degree less than $mp$: 
$$G_m=\left\{h(x)=\sum_{i=0}^{mp-1}a_ix^i\in \mathbf{k}[x]/(x^{mp}):a_{(i+1)p-1}=a_i^p, 0\leq i\leq m-1\right\}$$

Now we can state our result.

\begin{theorem}\label{thm:Regular}
    Let $M$ be an $A_{1,p}(R)$-lattice of rank $pm$, such that the reduction of $x$ modulo $p$ is regular nilpotent. 
    Then $M$ is isomorphic to $M_{f, g}$, for some $f, g$. Moreover, $M_{f_1, g_1}\cong M_{f_2, g_2}$ if and only if 
    $f_1=f_2, g_1=g_2\mod p$ and 
    $\frac{g_1-g_2}{p}\in G_m$.

\end{theorem}

\begin{proof}
    
Let $M$ be an $A_{1,p}(R)$-lattice of rank $mp$, such that the reduction of $x$ modulo $p$ is regular nilpotent. So $M/pM\cong\mathbf{k}[x]/(x^{pm})$ as $\mathbf{k}[x]$-modules. By Nakayama's lemma, $M$ is a cyclic $R[x]$-module. Thus, for some $f(x)$ of degree at most $mp-1$, we have $M\cong R[x]/(x^{mp}-pf(x))$ as $R[x]$-modules. 

Put $y\cdot 1=g(x)$,  where $\deg(g)\leq mp-1$.  Then it is immediate that $y=p\partial+g(x)$.  Hence $M\cong M_{f,g}$.

It remains to classify modules $M_{f, g}$ up to $A_{1,p}(R)$-isomorphisms. First, since the annihilator of $M_{f,g}$  
as an $R[x]$-module
is the ideal $(x^{mp}-pf)$, it follows that $f(x)$  is uniquely determined by the isomorphism class of $M_{f,g}$. 
It is straightforward to see that $M_{f,g_1}\cong M_{f,g_2}$  if and only if there exists $\phi\in\big(R[x]/(x^{mp}-pf(x))\big)^\times$ 
such that $g_1, g_2$ are related by the $p$-gauge transform:
$$g_1-g_2=p\phi'(x)\phi(x)^{-1}.$$

In particular, we see that $g_1=g_2 \mod p$ if such an isomorphism exists. Moreover, we see that $M_{f,g_1}\cong M_{f,g_2}$
if and only if $M_{f, 0}\cong M_{f, g_2-g_1}$. 

Next, suppose that $g\equiv 0\pmod p$, and let
$h\in\mathbf{k}[x]/(x^{mp})$ be the element corresponding to $g$ under
the natural identification
\[
p\big(R[x]/(x^{mp}-pf(x))\big)\cong\mathbf{k}[x]/(x^{mp}).
\]
Thus, we may write $h=g/p$ under this identification. The $p$-gauge
transformation condition implies that $M_{f,g}\cong M_{f,0}$ if and only if
\[
N_h=\big(\mathbf{k}[x]/(x^{mp}),\partial+h(x)\big)
\cong\big(\mathbf{k}[x]/(x^{mp}),\partial\big)=N_0
\]
as $A_1(\mathbf{k})$-modules.

Finally, a necessary condition for an isomorphism between $N_h$ and $N_0$
is equality of their $p$-curvatures: the actions by $y^p$ must match.
Since $y^p$ annihilates $N_0$, we must have $(\partial+h)^pN_h=0$.
Using the following identity in $A_1(\mathbf{k})$ (see \cite[Theorem~1.3]{B}), we have
$$(\partial+h(x))^p=\partial^p+h(x)^p+\partial^{p-1}(h(x)),\quad h(x)\in\mathbf{k}[x].$$
Writing $h=\sum_{j=0}^{mp-1}a_jx^j$, the coefficient of $x^{ip}$ in
$h^p+\partial^{p-1}(h)$ is $a_i^p-a_{(i+1)p-1}$ for
$0\le i\le m-1$, and all other coefficients vanish. Thus $h\in G_m$.

Conversely, if $y^pN_h=0$, then by Cartier descent (mentioned above)
there exist $\mathbf{k}[x^p]$-modules $L_h,L_0$ such that
$$N_h\cong \mathbf{k}[x]\otimes_{\mathbf{k}[x^p]}L_h,\quad
N_0\cong \mathbf{k}[x]\otimes_{\mathbf{k}[x^p]}L_0.$$
Since as $\mathbf{k}[x]$-modules both $N_h,N_0$ are isomorphic to a cyclic
module $\mathbf{k}[x]/(x^{mp})$, it follows immediately that
$L_h,L_0$ are isomorphic to $\mathbf{k}[x^p]/(x^{mp})$.
Hence $N_h\cong N_0$ as $A_1(\mathbf{k})$-modules. This completes the proof.

\end{proof}

\end{document}